\documentclass[12pt]{article}
\usepackage{amsmath}

\usepackage{times}
\usepackage{amsfonts}
\usepackage{amssymb}
\usepackage[english]{babel}

\newtheorem{theorem}{Theorem}[section]
\newtheorem{lemma}[theorem]{Lemma}
 \newtheorem{remark}[theorem]{Remark}

  \newtheorem{definition}[theorem]{Definition}
\newtheorem{proposition}[theorem]{Proposition}
\newtheorem{corollary}[theorem]{Corollary}

\newenvironment{proof}[1][Proof]{\textit{#1.} }{\hfill   \rule{0.5em}{0.5em}}
\def \qed { $\hfill \square$ \vskip 8 pt}

\def \endp { $\hfill \square$ \vskip 8 pt}

\renewcommand{\theequation}{\arabic{section}.\arabic{equation}}

\renewcommand{\rho}{\varrho}
\renewcommand{\epsilon}{\varepsilon}
\renewcommand{\phi}{\varphi}
\renewcommand{\theta}{\vartheta}
\newcommand{\e}{\epsilon}
\newcommand{\m}{\mu}

\newcommand{\grad}{\mathrm{grad}}
\newcommand{\s}{\sigma}
\renewcommand{\o}{\omega}

\newcommand{\la}{\lambda}

\renewcommand{\k}{\kappa}
\renewcommand{\d}{\delta}
\newcommand{\z}{\zeta}
\renewcommand{\a}{\alpha}

\newcommand{\n}{\nu}
\newcommand{\g}{\gamma}

\newcommand{\p}{\partial}

\newcommand{\Ric}{{\mathrm{Ric}}}

\newcommand{\wt}{\widetilde}

\renewcommand{\l}{{\lambda }}
\newcommand{\R}{{\mathbb{R}}}

\newcommand{\G}{{\Gamma}}

\newcommand{\M}{{\mathsf{M}}}

\renewcommand{\H}{{\mathsf{H}}}

\renewcommand{\L}{{\mathcal{L}}}

\renewcommand{\S}{\mathsf{S}}

\newcommand{\D}{{\nabla}}
\newcommand{\U}{{\mathcal{U}}}

\newcommand{\barint}
{\rule[.036in]{.12in}{.009in}\kern-.16in \displaystyle\int}

\renewcommand{\r}{\varrho}
\newcommand{\wh}{\widehat}

\date{ }

\begin{document}
\title{Liouville theorem, conformally invariant cones and  umbilical surfaces for
 Grushin-type
   metrics \thanks{MSC: 30C35, 53A30. Keywords: Riemannian 
geometry, Conformal maps, Liouville theorem, Umbilical surfaces.}}

\author{Daniele Morbidelli    }

\maketitle

\begin{abstract} We prove a classification theorem for conformal maps with respect to the control distance generated by a system of diagonal  vector fields in $\R^n$.
It turns out that    in many cases all such maps can be obtained as compositions  of    suitable dilations, inversions and isometries. Our methods involve a study of the singular Riemannian metric associated with the vector fields. In particular we  identify some conformally invariant cones related to   the Weyl tensor.  The knowledge of such cones enables us to classify all umbilical hypersurfaces.
 \end{abstract}


\section{Introduction} 
The principal purpose of this paper is to  classify  maps which are conformal with respect to
  the control  (Carnot--Carath\'eodory) distance  $d$
  generated by  a system of diagonal vector fields.
  Our principal result is that 
all such maps  
      are  compositions  of a restricted class of elementary conformal maps: 
isometries,  
   suitable dilations and  inversions 
    naturally associated with the distance $d$.
The form of these elementary maps will be explicitly identified. 


Consider in    $\M : = \R^p\times\R^q $  the diagonal vector fields
\begin{equation}
\label{gruco}
X_j=\frac{\p}{ \p{x_j}},\quad Y_\l = (\a+1)|x|^{\a} \frac{\p}{\p{y_\l}}, \qquad j=1,\dots, p,\quad \l=1,\dots q.
\end{equation}
Here $\a>0$ is  fixed.  Vector fields of the form
 \eqref{gruco} are usually referred to as Grushin vector fields and they are a 
subclass of the  diagonal vector 
fields studied by Franchi and Lanconelli  in \cite{FL}.  Denote by $d:\M\times\M\to\left[0,+\infty\right[$ the control
distance associated with
  the vector fields  in \eqref{gruco} 
 (see Subsection \ref{lasec},  or   
\cite{FL}  for a complete account).
We take here the following metric definition of conformal map. Let
$\Omega,\Omega'\subset \M$  be open sets. A homeomorphism
$f:\Omega\to \Omega'$ is  conformal  with respect to the metric
$d$ if there is a    function $u:\Omega\to\left] 0, +\infty\right[$ such that
\begin{equation}
\label{3051}
 \lim_{\z \to z} \frac{d(f(\z), f(z))}{d(\z, z)}= u(z)^{-1},
\end{equation}
for any $z=(x,y)\in \Omega$.   We say that $u$ is the
\emph{conformal factor} of $f$.

It is not difficult to check that the following maps are
conformal:
  \begin{equation}
\label{oooo}    (x,y)     \mapsto\Gamma(x,y)= (Ax, By+b ),\quad A\in O(p), \,\,B\in O(q),\,\, b\in \R^q;
\end{equation}
\begin{equation}\label{oooo1}
(x,y)   \mapsto
\delta_t(x,y):= (tx, t^{\a+1}y),\qquad  t>0.
\end{equation}
Maps of the form \eqref{oooo} are isometries. As the 
form of the vector fields $Y_\l$ suggests, 
 no translations in the variable 
  $x$ are admitted in \eqref{oooo}.   Note also 
  that  all the
   vector fields $X_j, Y_\l$ are homogeneous of degree 1 with respect to the anisotropic dilations   \eqref{oooo1}.

A less trivial
example of conformal map, which makes the model studied here
quite rich, is given by the following inversion. 
 Define the  ``homogeneous norm'' $
  \|z\|
   = \|(x,y)\|= \big(|x|^{2(\a+1)}+|y|^2\big)^{1/(2(\a+1))}$.  
   Then,
   for any $z \in \M\setminus\{(0,0)\} $, let
   \begin{equation}
\label{oooo2} \Phi(z) = \delta_{\|z\|^{-2}}z.
\end{equation}
The map \eqref{oooo2} is a  reflection in the homogeneous sphere
of equation $\|z\|=1$. It 
generalizes to to the present 
  setting
   the classical M\"obius
    inversion $t\mapsto |t|^{-2}t$, where $t $ belongs to an Euclidean space. See 
the discussion in
     Subsection  \ref{lasec}.
The conformality of the  map $\Phi$ was already recognized in \cite{MM}   by R. Monti and the author.

Compositions of the elementary  
maps described above provide easily more examples of conformal maps. 
 Our  main result states that, if $p\ge 3$, there are no further examples.
   \begin{theorem}\label{principale2}
Assume that  $p\ge 3$ and  $q\ge 1$.  Let   $\Omega,\Omega'\subset \M$ be  connected open sets. 
Let $f: \Omega\to \Omega' $
 be a  conformal homeomorphism in the metric  sense \eqref{3051}.
 Then  $f$   has the form
\begin{equation}\label{2553}
f(z) = \Gamma \Big( \delta_{t \|(x, y-b) \|^{s}}(x, y-b)\Big),
\end{equation}
for all $z=(x,y)\in \Omega$. Here
 $\Gamma$ is an isometry of the form  \eqref{oooo}, $t>0$, $b\in \R^q$ and $s=0$ or $-2$.
\end{theorem}
We immediately observe  that   the theorem is false
for $p=2$, $q\ge 1$. This is   a consequence of the fact
that the Riemannian
 metric $\wh g$ (see \eqref{cuc} below) is conformally
flat if $p=2$. See  Subsection \ref{noco} and  Remark
\ref{noleggio}. Case $p=1$ has been discussed by Payne \cite{PPP}.

The proof of Theorem \ref{principale2} requires Riemannian
arguments, because the control distance of the Grushin
  vector
fields is Riemannian away from the  (somehow small) set where $x$ vanishes.
Indeed,
 consider in $\M_0:= (\R^p\setminus\{0\})\times\R^q$, the  metric
 \begin{equation}\label{cuc}
 \wh g = |dx|^2 + \frac{|dy|^2}{(\a+1)^2|x|^{2\a}}.
 \end{equation}
   The vector fields introduced in \eqref{gruco} form an orthonormal
   frame in $(\M_0, \wh g)$.
   It is easy to realize
    that their  control distance agrees with the
Riemannian distance $d_{\wh g}$ associated with $\wh g$ (lengths
of curves are the same).  Moreover, by a result of
Ferrand \cite{F},  a conformal homeomorphism in a smooth Riemannian manifold  must be  smooth. 
Thus, given a homeomorphism
  $f:\Omega\to \Omega'$,    where $\Omega,\Omega'\subset \M$, if $f$ satisfies  
\eqref{3051}, then   it is smooth in $\Omega\cap\M_0$ and it 
 must satisfy    the
Cauchy--Riemann system
\begin{equation}
\wh g(f_*U, f_*V) = u(z)^{-2}\wh g(U, V),
\label{crr}
\end{equation}
 for all vector fields   $U,V$ supported   $\Omega\cap\M_0$ and for a suitable conformal factor $u$.

In view of the discussion above, it turns out that 
our main result Theorem \ref{principale2}
follows immediately from the following Liouville theorem for  the manifold    $(\M_0, \wh g)$. 
   \begin{theorem}\label{principale}
Let  $p\ge 3$ and  $q\ge 1$.  If   $\Omega,\Omega'\subset \M_0$
 are  open and
connected and  $f:\Omega\to\Omega' $ is a smooth
 diffeomorphism satisfying the Cauchy--Riemann
  system
 \eqref{crr}, then
   $f$   has the form \eqref{2553}.
\end{theorem}

In order to prove Theorem \ref{principale},
 we use the metric $g= (\a+1)^2|x|^{2\a}\wh g$, which belongs to
  the same conformal class of $\wh g$ and makes computations easier.
A standard way to study conformal maps on Riemannian manifolds
 starts from the interpretation of the  transformation formula for
  the Ricci tensor under  conformal
  changes of metrics as a tool to obtain a  system of partial differential 
equations for 
the conformal factor $u$. Indeed,   given a   metric $g$, letting $\wt g = u^{-2}g$, then we
  have the classical   formula 
 \begin{equation}\label{uuvv}
 \Ric_{\wt g}  = \Ric_g  + (n-2)u^{-1}\D^2 u -  u^{-2}\big\{
(n-1)|\D u|^2 - u\Delta u \big\}g.
\end{equation}
Here $\D^2 u$ denotes the Hessian in the Levi Civita connection
$\D$ of the metric $g$, while $|\D u|$ is the length of the gradient and
$\Delta u$ the Laplacian. Next, if $f :\Omega\to \Omega'$
is a
conformal diffeomorphism
in $(\M_0,g)$,  which means, by definition,  $f^*g = u^{-2}g$ for some function $u$, then $u$ satisfies
\begin{equation}\label{pazziata3}
\begin{split}
&
\mathrm{Ric}(f_*U, f_*V)
    \\&
= \mathrm{Ric}(U, V)+ (n-2)u^{-1}\D^2  u(U, V) - g(U, V)u^{-2}\big\{ (n-1)|\D u|^2 -
u\Delta u \big\},
\end{split}
\end{equation}
for any  pair of  vector fields $U, V$ in $\Omega$. Here $\Ric = \Ric_g$.
We will show that, if $p\ge 3$,
 then all  solutions $u$ of   \eqref{pazziata3}
 have the same form as the conformal factor of a suitable composition of maps of 
the form 
 \eqref{oooo}, \eqref{oooo1} and \eqref{oooo2}. This  will reduce the proof to a 
classification of
 local  isometries of $g$, which is given in Section \ref{duedue}.

This strategy
 can be easily pursued in  the Euclidean case and it
 reduces to a few lines 
thanks to the Ricci flatness of the Euclidean space.  See \cite[Chapter 6]{SY}. (See also \cite{Ber} or the 
 Liouville's  paper \cite{Liou}, for analytical proofs 
   based on differentiation of the Cauchy--Riemann system.) 
        Indeed, in the Euclidean case,
   Ricci flatness makes  system \eqref{pazziata3} an easy
   overdetermined system in the only unknown  $u$, whose
   solutions are particular quadratic
   polynomials.  In our case, the metric
   $g$ is not   Ricci flat.  Therefore it is not clear how to manage
 the Ricci terms in
   \eqref{pazziata3}, especially the one on the left--hand side.

In order to make viable  system \eqref{pazziata3}, we introduce  the
following  conformally  invariant
 cones $\U_P\subset T_P \M_0$.
 Assume that   the dimension of $ \M_0$ is greater than 3. 
 Let $R, \Ric$ and $ \mathrm{Scal}$ be the
 Riemann, Ricci and scalar curvature of $g$, respectively. Let
\begin{equation} \label{wwi}W = R +\frac{1}{n-2}(\Ric\odot g)
-\frac{\mathrm{Scal}}{2(n-1)(n-2)}(g\odot g)
\end{equation}
be the Weyl tensor.  Here $(h\odot s)_{abcd}:=
h_{ad}s_{bc}+h_{bc}s_{ad}-h_{ac}s_{bd}-h_{bd}s_{ac}$ denotes the
Kulkarni-- Nomizu product of symmetric 2-tensors.
 Then
define
 \begin{equation*}
\begin{split}
   \U_P = \{  & X\in T_P \M_0 : W(X, Y, U, V ) =0\text{ for all $Y,U, V  $}
\\ & \text{such that   $X,Y,U, V  $ are pairwise orthogonal with respect to $g$   }\}.
\end{split}
\end{equation*}
 The conformal invariance
of the Weyl tensor (if $\wt g= u^{-2}g$, then $ W_{\wt g} =u^{-2}
W_g $)  trivially implies that if $f:\Omega\to
\Omega'$ is a conformal diffeomorphism between subsets of $
\M_0$, then
\begin{equation*}
f_* (\U_P) = \U_{f(P)}, \qquad \forall\,\, P\in \Omega.
\end{equation*}
In the case of the Grushin metric, the cones $\U_P$ will be
determined explicitly  at any $P\in \M_0$. They have a very clear structure in suitable cylindrical coordinates. 
Observe that the invariants   $\U_P$ may  also   be used to study conformal maps
in different  Riemannian manifolds, provided the Weyl tensor has a non trivial structure. See, for example,  Remark \ref{cito2}. 
 We also
 mention that  
 different
conformally invariant subsets (actually subspaces) of the tangent
space constructed from the Weyl tensor were used by  Listing 
\cite{L}.

 With the explicit form of  the cones $\U_P$ in hand, it becomes possible
  to deal with Ricci terms in   system
 \eqref{pazziata3} and ultimately to solve it.
  Then Theorem \ref{principale}
   follows in a rather standard way. It turns out that
 all conformal maps  preserve the Ricci tensor,  in the sense that
   $\mathrm{Ric}_g(f_*U, f_*V)
= \mathrm{Ric}_g(U, V)$, for all vector fields $U, V$. Therefore
conformal maps
 are    Liouville maps in the
language  of \cite{KR}
 and, in particular,
   M\"obius maps in the terminology
    of \cite{OS}.  Here the choice of the metric $g$ in the  conformal
 class $[\wh g]$ is important.

Although our method of classification   could be less efficient than other techniques, like the study of 
conformal Killing vector fields (see Payne  \cite{PPP}, in the case  $p=1$), it should be emphasized
that our approach provides more.
Indeed, the   study of 
the cones $\U_P$   is inspired by  the  
fact that
conformal maps   must    send    umbilical
hypersurfaces  to umbilical hypersurfaces,  see
 e.g. 
  \cite{KP}. 
 Indeed, given a manifold  $(\M_0,g)$, 
of dimension  at least $ 4$, a standard obstruction
 to the existence of an umbilical surface $\Sigma$ with
given normal $N\in T_P \M_0$ at a point $P\in \Sigma$
is provided by Codazzi equations (see \cite{ON}),
which for an umbilical hypersurface of curvature $\k$ with respect to a unit normal $N$ become
  \begin{equation}
  \label{paippo}
 g(V, Z) ( U \k)   -  g(U, Z) (V \k)    = R(N, Z, U, V),
\end{equation}
for all $U, V, Z$ tangent to $\Sigma$.   If  $U, V$ and $ Z$ are 
pairwise
orthogonal,  then
 $R(N, U, V, Z)=0$. Moreover,   $R(N, U,
V, Z )= W(N,U, V, Z)$, by the form \eqref{wwi} of Weyl tensor. Therefore,
equation \eqref{paippo} shows that a normal vector at $P$
 to an umbilical
surface cannot belong to $\U_P$. 

After the cones $\U_P$  are known, we are able to classify in Section \ref{ooon}  all umbilical hypersurfaces in    $(\M_0,g)$ for
$p\ge 3$. It turns out that they are rather rare,  
while for $p=2$ the situation is different, see  Remark \ref{ombrello2}.   
Here is our result.
\begin{theorem}
\label{ombelo}
Let  $p\ge 3$, $q\ge 1$.  Then any connected umbilical
hypersurface in  $(\M_0,g)$ is contained in one among the following:

\begin{description} \item   [{\rm(A1)}]
a homogeneous sphere of equation
$
|x|^{2(\a+1)} + |y-b|^2 = c^2, \quad b\in \R^q, c>0;
$

\item   [{\rm(A2)}]   a plane of equation
 $
\langle  a, y\rangle=c,
 $ where
 $a\in \R^q$, $c\in \R$;

 \item   [{ \rm  (B)}] a plane of equation
$
\langle a,x\rangle=0,
 $ where $a\in \R^p$.
 \end{description}
 \end{theorem}

Observe  that the choice of the metric $g$ in the conformal class $[g]$ ensures that  all umbilical
surfaces have constant curvature, see Remark \ref{ombrello2}.

  Grushin--type geometries  have been recently
    rather studied. They pose interesting problems
   from the point of view of nonlinear analysis, sharp inequalities
    and search for symmetries related to the degenerate elliptic operator
\begin{equation}
\label{numero2}
\Delta_\a:=  \Delta_x  + (\a+1)^2 |x|^{2\a}\Delta_y .
\end{equation}
 See, for example,  the papers
   \cite{B,M,T,GV,BG,LP1,LP2}, just to quote a few.
   The conformal inversion $\Phi$ in \eqref{oooo2} is used in  \cite{MM},  
      in order to construct
   a Kelvin--type transform for a semilinear equation with critical nonlinearity
    of the form
 $-\Delta_\a u = u^{r}$, for a suitable $r>1$.
Our motivation for a better understanding of these
 conformal maps stems from the mentioned paper.



Concerning  Liouville--type theorems in sub--Riemannian geometry,
    we mention the seminal papers by Kor\'anyi and Reimann \cite{KR1,KR2},
   where conformal maps in the Heisenberg group were classified (see also 
\cite{C,T}). More
 rigidity
    results  are contained in    \cite{P,RR,R}.  
The quoted paper are   in the setting of Lie groups.

The  plan of the paper is the following. In  Section \ref{duedue},
we discuss some preliminary facts.  The metric $g$,   its
isometries (Subsection \ref{noco}),  and the  conformal inversion
$\Phi$ (Subsection \ref{lasec}).
 In    Section \ref{tretre} 
 we
 prove the  Liouville theorem. 
We first study the cones $\U_P$,  in Subsection \ref{trio2};
 then, in Subsection \ref{2952}, we
 solve system 
\eqref{pazziata3};   finally we show, in Subsection \ref{triotre}, how the proof 
can be quickly concluded in view of the explicit knowledge 
of isometries.    Section \ref{ooon}  is devoted to the  classification of umbilical surfaces.
Finally, we included a short appendix with some standard formulas on warped metrics in Riemannian manifolds.

\medskip\noindent\bf Notation. \rm  Given a Riemannian metric $g$, we
denote by $\D$ the associated
  Levi Civita connection and by  $R(X, Y)U=
\D_X\D_YU-\D_Y\D_XU-\D_{[X, Y]}U$ the curvature operator. We let
$R(U, V,X,Y )=g(U, R(X, Y)V)$, so that $\Ric(X, Y)= \mathrm{trace}
\{V\mapsto R(V, Y)X \}$. Moreover, $\langle\cdot,\cdot \rangle$
denotes Euclidean scalar product. Surfaces have codimension 1 and
are orientable and connected. Unless otherwise stated, Latin
indices $i,j,k$ run from $1$ to $p$, while $\l,\m,\s$ go from
$1$ to $q$. For typographical reasons, we write
$\p_j$ or $\p_{x_j}$ instead of  $\frac{\p}{\p x_j}$ and we use
the an analogous notation for   $\frac{\p}{\p y_\l}$.
 Summation with
respect to a repeated index (in the pertinent range) is sometimes
omitted.

\medskip\noindent\bf
Acknowledgements. \rm 
     It is a pleasure  to thank
Roberto Monti: the problem and ideas contained in this paper
originated from  our discussions and our shared   research
activity. In particular,  I am indebted to him for his help in
the discussion of  system \eqref{pazziata3}.

 I also thank Zolt\'an Balogh and Luca Capogna who 
kindly indicated  same references to me.

 \section{Preliminary facts on  Grushin geometry}
 \label{duedue}\setcounter{equation}{0}
\subsection{The Grushin metric} \label{noco}  In the  conformal class of  $\wh g$ we  choose the following metric  $g$ in $\M_0$:
     \begin{equation}
\label{straccio}
\begin{split}
  g & = (\a+1)^2|x|^{2\a}|dx |^2 +  |dy|^2.
\end{split}
\end{equation}
It turns out that    $g$ is better than  $\wh g$  for our purposes. Observe that, if $p=2$, then the local holomorphic change of variable $x^{\a+1}=\xi$, or alternatively  formula \eqref{uuvv} show that, for any $\a\ge 0$, the metric
$(\a+1)^2|x|^{2\a}|dx|^2$ is flat if  $x\in\R^2\setminus\{  (0,0)\}$.
 Therefore, $g$
is flat for $p=2$, $q\ge 1$.

Let $\r=|x|$, $\theta= \frac{x}{|x|}\in \mathbb{S}^{p-1}\subset \R^p$.
 Then, $(\a+1)^2|x|^{2\a}|dx|^2= (\a+1)^2\{
\r^{2\a}d\r^2 + \r^{2(\a+1)}|d\theta|^2 \}$, where
 $|d\theta|^2$ is the standard
metric on the sphere $\mathbb{S}^{p-1}$.
Moreover, letting $\r^{\a+1}= r $,
 it follows quickly that 
\begin{equation}
\label{harpo4} g=  dr^2 + |dy|^2+ (\a+1)^2   r^2 |d\theta|^2.
\end{equation}
  Using 
the notation
\begin{equation*}
\begin{split}
\H & =  \{ (r,y) \in \left]0,+\infty\right[\times \R^q\} , \qquad g_\H = dr^2 + |dy|^2,
\\
\S &= \mathbb{S}^{p-1}=\{ \theta\in \R^p:|\theta|=1\}, \text {  $g_\S=|d\theta|^2$},\end{split}
\end{equation*}
we can write 
the  manifold $(\M_0, g)$
   as a warped 
product $\H\times_w \S$, with warping function $w(r,y)= (\a+1)r$. Briefly, $g=g_\H+w^2 g_\S$. See the appendix for some standard facts about warped products. See \cite[Chapter 7]{ON} for a complete   introduction.

For any $P=(r,y,\theta)\in \H\times \S$, decompose any $U\in T_P (\H\times\S) $ as
\begin{equation}
\label{aeio}
  U= U_\H + U_\S\in T_P\H \oplus T_P\S,
\end{equation}
where $T_P\S$ and $T_P\H$ denote the   lifts 
 at $P$ of the tangent spaces $T_{\theta}\S $ and  $T_{(r,y)}\H$, respectively.

 Next we describe the connection in the warped model.
  Denote
by  $\D^\H$,  $\D^\S$ and  $\D$  the Levi Civita connections on $\H, \S$ and $\H\times_w\S $, respectively.
 Since the  factor $\H$ is
Euclidean, by \eqref{alfa},
  covariant derivatives are Euclidean,
namely (in  the notation  $\p_r = \p/\p r$ and $\p_\l= \p/\p{y_\l}$)
  \begin{equation}
  \label{cane}
  \D_{\p_r}\p_r
 =\D_{\p_r}\p_{\l}= \D_{\p_{\l}}\p_{\m}=0,\qquad \l,\m=1,\dots,q.
 \end{equation}
  Then $\D^2 u(\p_r, \p_r) =
\p_r^2u$, $\D^2 u(\p_r, \p_{\l})=\p_r\p_{\l} u$  and $\D^2 u(\p_{\l},
\p_{\m})=\p_{\m}\p_{\l} u $.
Moreover, again by \eqref{alfa}, given $X$,  lifting of a vector field on the sphere, then
\begin{equation} \label{fiorfiore}
\D_{\p_r} X = \frac{ 1   }{r}  X \quad \text{and}\quad \D_{\p_{\l}} X =0,\quad  \l=1,\dots, q.
\end{equation}
Let
 $u = u(r,y,\theta)$ be a scalar function. Then
$
|\D u|^2 =   (\p_ru)^2  + |\D_y u|^2+ (\a+1)^{-2} r^{-2}|\D^\S u|^2.
$ Moreover,
given $X,X'$ on the sphere,  \eqref{beta} provides 
\begin{equation}\label{sforo}
\D^2 u(X, X')  =  (\D^\S)^2 u (X, X')+ \frac{u_r}{r}g(X, X')\qquad\text{and}
\end{equation}
\begin{equation}\label{tato}
\begin{split}
\Delta u  &
 =   u_{rr}+ \frac{p-1}{r}u_r+\Delta_y u   +\frac{1}{
 (\a+1)^2 r^2}\Delta_\S u.
\end{split}
\end{equation}

In order to compute the curvature,
 note that in our case  $\D^2 w=0$ ($w$ is linear in $r,y$).
  Therefore only the third line in
\eqref{beta3} gives nonzero terms. A short computation
using the curvature of the
standard sphere $\S=\mathbb{S}^{p-1}$, $R^\S (V_1, V_2, V_3, V_4)=
g_\S (V_1,V_3 ) g_\S (V_2, V_4)- g_\S (V_1, V_4) g_\S (V_2, V_3)$,
gives
\begin{equation}
\label{mimmi}
\begin{split}
 R(U, &  V, X,  Y) = R(U_\S, V_\S, X_\S, Y_\S)
 \\& =  -\a(\a+2)(\a+1)^{-2} r^{-2} \Big\{ g  (U_\S,X_\S) g  (V_\S, Y_\S)- g
(U_\S, Y_\S) g (V_\S, X_\S)\Big\},
\end{split}
\end{equation}
where all the vectors    $U, V, X, Y$  are decomposed as in
\eqref{aeio}. We see again that
the manifold is flat for $p=2$ ($
\mathbb{S}^{p-1}= \mathbb{S}^1$ and the curly bracket in
\eqref{mimmi} vanishes).   Contracting,
\begin{equation}
\begin{split}
\label{riciao}
\mathrm{Ric} (U, V)& = \mathrm{Ric} (U_\S, V_\S) = -\a(\a+2)(p-2)
(\a+1)^{-2} r^{-2} g (U_\S, V_\S) \quad\text{and}
\\ \mathrm{Scal} & =-\a(\a+2)(p-2)(p-1)(\a+1)^{-2} r^{-2} .
\end{split}
\end{equation}

Next we
classify all local isometries of $\H\times_w\S$ for $p\ge 3$.
\begin{proposition}\label{ioio}
Let $p\ge 3$. Let $\Omega\subset \M_0$ be a connected open set. Let $f:\Omega\to \Omega'\subset\M_0$ be a   local isometry in the metric $g$.
Then $f$ is a restriction of a map of the   form 
\[
(x,y)    \mapsto  (Ax, By+b ),\qquad 
\]
where $A\in O(p),$ $B\in O(q)$ and $b\in \R^q$.
\end{proposition}
\noindent\it Proof. \rm Write the map as  $(x,y)\mapsto (\wt x (x,y),
\wt y(x,y))$. Since isometries preserve scalar curvature,
 \eqref{riciao} gives  $|\wt x(x,y)| =
|x|$ for all $(x,y)\in\Omega$.  Here the choice $p\ge 3$ is
crucial.

Introduce the notation $\Sigma_\r= \{(x,y):|x|=\r\}$. Next we claim that, for any  $\r>0$,  the restriction of the map $f$ to the set $\Omega\cap\Sigma_\r$ (provided the latter is nonempty) has the form
\begin{equation}
\label{fui}
(x,y)\mapsto (A(\r) x, B(\r) y + b(\r)),
\end{equation}
where $A(\r)$, $B(\r)$ are orthogonal and  $b(\r)\in \R^q$. We assume without loss of generality that $\Omega$ is a product of the form $\Omega = \{ \r\in (\r_0,\r_1), |y-y_0|<\e_0, |\theta -\theta_0|<\e_1 \}$, so that $\Omega  \cap\Sigma_\r $ is connected. 
Observe that the metric on $\Sigma_\r$ has the form
 $(\a+1)^2 \r^{2(\a+1)}|d\theta|^2  + |dy|^2$, the product of a sphere (of dimension at least $  2$, because $p\ge 3$) with an Euclidean space.
 Therefore the claim follows  from the standard fact than  a local isometry of a product of space forms of different curvature must be a product map of isometries of the factors (this fact can be easily proved
  by means of isometric invariance of sectional curvatures).

Finally we prove that $A, B, b$ are constant in  $\r$. Take a point $z=(x,y)$. The normal vector $(\p_{\r})_{z}$ at the point $z=(x,y)$
 to the surface   $\Sigma_\r $ 
 is sent by
$f_*$ to the vector  $\pm(\p_{\r})_{\wt z}$, normal to the same surface at the point  $f(z)= \wt z= (\wt x, \wt y)$.  Since  scalar curvature is increasing as $\r$ increases, the sign must be $+$. Moreover, by the chain rule, 
\[
f_*((\p_\r)_z)
 = \p_\r \big(A_j^k x^j\big) (\p_{k})_{\wt z} + \p_\r
\big( B_\m^\s y^\m +b^\s \big) (\p_{\s})_{\wt z},
\]
where we sum for $j,k=1,\dots, p$ and for   $\s,\m=1,\dots,q$.
Therefore, the second term, the one with derivatives in $y$, must be   zero.
Thus, $(\p_\r B )y + \p_\r b=0$. Differentiating in  $y$,  we get $\p_\r B =0$. Then  $\p_\r b=0$. We have proved that  $B$ and  $b$ are constant.

Finally, we look at the first term. Recall that $\r=|x|$, so that  $\p_\r x^j =  x^j/ \r$. Then,
\[
\begin{split}
 (\p_\r)_{\wt z} &=f_*(\p_\r)_z = \p_\r \big(A_j^k x^j\big)
(\p_k)_{\wt z}
 = \Big((\p_\r  A_j^k ) x^j +  A_j^k \frac{x^j}{\r} \Big)(\p_k)_{\wt z}
\\&
=\Big((\p_\r  A_j^k ) x^j + \frac 1\r \wt x^k \Big)(\p_k)_{\wt z}
 =
(\p_\r  A_j^k ) x^j (\p_k)_{\wt z}
+ (\p_\r)_{\wt z}.
\end{split}
\]
Thus,
$
(\p_\r A_j^k ) x^j = 0,
$
which gives (differentiate in $x$)  $\p_\r A_j^k =0$. The proof is concluded.
\endp

 \subsection{Control distance and conformality of  the inversion map. }
\label{lasec} In this subsection we show that  inversion is conformal. The same result has been proved in \cite{MM}, but here
 we provide a shorter proof, using  the warped model.
  Let $\Phi(z)= \d_{\|z\|^{-2}}z$. Our aim is to check that, for any $z\neq (0,0)$,
\begin{equation}\label{1853}
\lim_{\z\to z}\frac{d(\Phi(\z), \Phi( z))}{d(\z,  z)}= \|z\|^{-2}.
\end{equation}

 Before proving \eqref{1853},
 we briefly recall  the definition
 of
 control distance associated with the vector fields $X_j, Y_\l$, $j=1, \dots, p, 
\l=1, \dots, q$.
 See \cite{FL}, see also \cite{NSW}.
An absolutely continuous path $\gamma:[0,T]\to \M$ is admissible if it  satisfies almost everywhere  $\dot\g = \sum a_j X_j(\g) + \sum b_\l Y_\l(\g)$  for suitable measurable functions $a_j, b_\l:[0,T]\to \R$. Define, for $z, z'\in \M
$,
$
d(z, z')= \inf\int_0^T\sqrt{|a|^2 + |b|^2},
$
where the infimum is taken among all the functions $a_j, b_\l$
such that the corresponding path $\g$ is admissible and connects
$z$ and $z'$.

We   prove conformality by means of a suitable Cauchy-Riemann system.
Indeed
 we prove that
\begin{equation}\label{a2}
\wh g(\Phi_*U, \Phi_*U) = \|z\|^{-4}\wh g(U, U),
\end{equation}
for all $U\in \mathrm{span}\{ X_j, Y_\l : j=1,\dots, p, \l=1,
\dots, q\}$, $z=(x,y)\neq( 0,0)$.

We first prove  \eqref{a2}   in the set 
  $\M_0$, namely where $|x|>0$.
 In the warped model of Subsection \ref{noco},   the map $\Phi$ takes the form
 \[
\Phi(r,y,\theta) = \big ( \phi (r,y), \theta \big), \quad \text{ where }\quad
\phi(r,y) =  |(r,y)|^{-2} (r,y)\]
is an Euclidean M\"obius map.
The metric $g$ at $P:=(r,y,\theta )$ is $dr^2 + |dy|^2 + (\a+1)^2 r^2 |d\theta|^2$, while in $\Phi(P) = (\phi(r,y),\theta)$ it has the form
  $dr^2 + |dy|^2 + (\a+1)^2  |(r,y)|^{-4} r^2   |d\theta|^2$. Therefore, if
  we decompose, as in \eqref{aeio}, $U = U_{\H} + U_\S$, we get
 $
\Phi_*(U_\H + U_\S) = (\phi_* U_\H) + U_\S\in T_{\Phi(P) }\H\oplus
 T_{\Phi(P)}\S.
 $
 Thus, since $T\H$ and $T\S$ are orthogonal,
 \[
 g_{\Phi(P)} (\Phi_*U, \Phi_*U) =   g_{\Phi(P)} (\phi_*U_H, \phi_*U_H)  +
  g_{\Phi(P)}  ( U_\S, U_\S),
 \]
where, in order to be safe,  we used the  slightly 
cumbersome notation $g_{\Phi(P)}$
 to indicate the metric $g$ at the
point $\Phi(P)$. Next look at the first term. By the properties of
Euclidean M\"obius maps,  we have 
\begin{equation}
\label{didos}
 g (\phi_*{U_\H}, \phi_* U_\H) =  |(r,y)|^{-4}   g(U_H, U_\H)  .
\end{equation}
Moreover, looking at the second term, since the metric at the image point   $(\phi(r,y), \theta)$
is  $dr^2+ |dy|^2 + (\a+1)^2 |(r,y)|^{-4} r^2
g_\S$ (here $g_\S= |d\theta|^2$), we have
\begin{equation}\label{didas}
g_{\Phi(P)} (U_\S, U_\S)  = (\a+1)^2 |(r,y)|^{-4} r^2 (g_\S)_\theta(U_\S, U_\S)=|(r,y)|^{-4} g_{P}(U_\S, U_\S).
\end{equation}
Putting together the three formulas above,
\begin{equation}
\label{ciccia2}
g(\Phi_*U, \Phi_*U) = |(r,y)|^{-4}g(U, U) =
\|z\|^{-4(\a+1)}g(U, U) ,
\end{equation}
which will be referred to in Section \ref{tretre}. Since $ g= (\a+1)^2|x|^{2\a} \wh g = (\a+1)^2 r^{2\a/(\a+1)} \wh g $, we also get
\begin{equation}\label{a3}
\wh g(\Phi_*U, \Phi_*U) = \|z\|^{-4}\wh g(U, U),
\end{equation}
for all   vector field  $U$ in $\M_0$. Hence  \eqref{a2} is proved   at any point of $\M_0$.

Next  we prove \eqref{a2} at   points of the form  $(0,y)$, $y\neq 0$.
Here we may work in Cartesian coordinates.   Observe that
$\Phi(0,y)= (0, |y|^{-2} y)$ and
$\p_{x_j}(\|z\|)\Big|_{(0,y)}= 0$. Therefore it is easy to see that
$
\Phi_*(\p_{x_j})_{(0,y)} = |y|^{-2/(\a+1)}(\p_{x_j})_{(0, |y|^{-2} y)}.
$
Thus, if  $U= U^j(\p_{x_j})_{(0,y)}$,
\begin{equation}
\label{b}
\wh g (\Phi_* U,\Phi_* U)= |y|^{-4/(\a+1)} \wh g(U, U).
\end{equation}
Equations \eqref{a3} and \eqref{b}   together complete the proof of \eqref{a2}.

In order to prove conformality starting from \eqref{a2}, use the following routine argument. Take a point  $z_0\neq 0$. Let $z$
be a close point and denote $\e =d(z, z_0 )$. Take an arclength geodesic $\gamma:[0,\e]\to \M $,
$\gamma(0) =z_0 $, $\gamma(\e) = z $,  with $\dot \gamma(t) =
a_j (t)X_j(\gamma(t)) + b_\l(t) Y_\l(\gamma(t))$ and
$|a(t)|^2+|b(t)|^2=1$ at almost all $t$. We may assume that
$\gamma  $ does
 not touch $(0,0)$, provided $\e$ is small enough.
Then
\begin{equation}
\label{sss}
d\big( \Phi(z_0), \Phi(z)\big)\le
\int_0^\e\sqrt
{\wh g_{\Phi(\g)}(\Phi_*\dot\g, \Phi_*\dot\g)}dt = \int_0^\e \|\gamma\|^{-2} dt,
\end{equation}
because $\wh g (\dot\g, \dot\g)=1$ almost everywhere.
As $\e\to 0$ we get $\displaystyle{\limsup_{\e\to 0}}\frac{d\big( \Phi(z_0), \Phi(z)\big)}{d\big( z_0,z \big)}\le  \|z_0\|^{-2}$. The same argument applied to $\Phi^{-1}$ provides  equality \eqref{1853} at any point $z_0\neq (0,0)$.

\section{Proof of the Liouville theorem}\label{tretre}
 \setcounter{equation}{0}
In this section we first study the cones $\U_P$
 of the metric $g$. 
Then we use their form  to find all admissible conformal factors $u$ of a conformal map in $(\M_0, g)$, for $p\ge 3$. At the end of the section we give the easy argument which concludes
the proof of Theorem \ref{principale} and hence of Theorem \ref{principale2}.

\subsection{The cones   $\U_P$  for the   metric $g$.}\label{trio2}
In the following proposition we  identify  the cones $\U_P$ defined in the introduction. We use the warped metric \eqref{harpo4}.

\begin{proposition}\label{tasti}
Let   $p\ge 3$. Then, for any  $P\in \H\times_w \S$, we have
\[
\U_P = \{ X \in T_P(\H\times_w \S) : |X_\H|\,\, |X_\S|=0\} = T_P\S\cup T_P\H.
\]
\end{proposition}
Observe that, if    $p=2$ and $q\ge 2$, then we have  $\U_P = T_P\M_0$,
all the tangent space, because the metric is flat.

  Proposition  \ref{tasti},   
together with a continuity argument, immediately gives corollary below, whose  easy proof is omitted.
\begin{corollary}
Let $\Omega\subset \H\times\S$ be a connected open set.
Let 
$f:\Omega\to f(\Omega)\subset \H \times\S$ be a  conformal
diffeomorphism in the metric $g=g_\H + w^2 g_\S$. Assume that $p\ge
3$. Then, either
\begin{equation}\label{caso1}
\begin{cases}
    f_*(T_P\S)= T_{f(P)}\S
    \\
    f_*(T_P\H)= T_{f(P)}\H
\end{cases} \forall P\in \Omega,  \end{equation}
or
\begin{equation}\label{caso2}
\begin{cases}
    f_*(T_P\S)= T_{f(P)}\H
    \\
    f_*(T_P\H)= T_{f(P)}\S
\end{cases} \forall P\in \Omega. \end{equation}
Correspondingly, in cylindrical coordinates $(r,y,\theta)$, 
the map is a product of one between the following types:
\begin{equation}
\label{promo} (r,y,\theta)\mapsto (\wt r(r,y), \wt y (r,y),\wt
\theta(\theta)),
\end{equation}
or
\begin{equation}
\label{zione}
(r,y,\theta)\mapsto (\wt r(\theta), \wt y (\theta),\wt \theta(r,y)). \qquad
\end{equation}
\end{corollary}
Observe that, for dimensional reasons,  \eqref{caso2}  and the corresponding \eqref{zione} may happen only if   $ \S$ and  $\H$ have the same dimension, namely when
$p-1=q+1$.

\begin{remark}\rm
\label{noleggio} Since in case $p=2$ the metric $g$ is flat, it is easy to realize that in this situation there are conformal maps which  satisfy neither \eqref{caso1}, nor \eqref{caso2}. More precisely,
 given any point $P$ and any $X, Y\in T_P\M_0$, there is a local isometry $f$ 
around $P$ such that $f(P)=P$ and $f_*X=Y$.
\end{remark}

\noindent \it Proof of Proposition  \ref{tasti}.
  \rm Observe first
that  $W(X, U, V, Z )= R(X, U, V, Z)$, provided   $X, U, V, Z$
form an orthogonal family. This follows from \eqref{wwi}.

The proof will be accomplished in two steps. 

\smallskip

\it Step 1.  \rm If   $X= X_\H + X_\S\in T_P\H\oplus T_P\S$ with  $|X_\H|\neq 0 $ and $|X_\S|\neq 0$, then $X\notin \U_P$.

\it Step 2. \rm If  $X= X_\H\in T_P\H$ or $X= X_\S\in T_P\S$, then  $X\in \U_P$.

\smallskip {\it Proof of Step 1. }
Write $X= X_\H +X_\S\in T_P\H +T_P\S$. Recall that both $X_\H$ and $X_\S$ are nonzero. Take two nonzero vectors
$X^\bot_\S\in T_P\S$ with $g(X_\S, X^\bot_\S)=0$ and
$X^\bot_\H\in T_P\H$ with $g(X_\H, X^\bot_\H)=0$. This choice is 
possible, because $\dim T_P\S\ge 2$ ($p\ge 3$) and $\dim T_P\H=q+1\ge 2$. Then take
$
V = X_\S - c_1 X_\H,
$
 where $c_1$ is such that  $ g(X, V)=0$, and
$
U   = X_\S^\bot + X_\H^\bot,$ $  Z = X_\S^\bot - c_2 X_\H^\bot,
$ where $c_2$ is such that  $ g(U, Z)=0.$
Then $X,U,  V, Z$
 form an orthogonal family and moreover, by \eqref{mimmi},
$
R(X, U, V, Z) = R(X_\S, X_\S^\bot,X_\S, X_\S^\bot)\neq 0.
$
Step 1 is proved.

\medskip
\it  Proof of Step 2.  \rm If  $X = X_\H\in T_P\H$ and   $X, U, V,
Z$ form an orthogonal family, then   $W(X, U, V, Z)=R(X, U, V,
Z)=0$, by \eqref{mimmi}.

If   $X= X_\S\in T_P\S$,
 take  $U, V, Z$ orthogonal triple, where all the vectors $U, V $ and 
 $ Z$ are  orthogonal to
  $X$.
There are two cases.

 First case: all the vectors
   $U, V, Z$ have nonzero projection along $T_P\S$, $U= U_\H + U_\S $, $V= V_\H + 
V_\S $ and  $Z= Z_\H+Z_\S$, with
  $|U_\S|\,|V_\S|\,|Z_\S|\neq 0$. But then, since $U, V, Z, X$ are orthogonal and 
$X_\H=0$,
 all $U_\S, V_\S$ and $ Z_\S$ must be orthogonal to $X_\S$. Hence, by \eqref{mimmi},
 \[
 \begin{split}
 W(X,U, V, Z)& = R(X, U, V, Z)=  R(X_\S, U_\S, V_\S, Z_\S)
 \\& =-\a(\a+2)(\a+1)^2 r^2 R^\S
 (X_\S, U_\S, V_\S, Z_\S)
  =0,
 \end{split}
 \]
 by elementary properties of the curvature $R^\S$ of the sphere.

 Second case:
 at least one among the vectors $U, V, W$ has zero projection along   $T_P\S$. Then $W(X, U, V, W) = R(X, U, V, W)=0$, by \eqref{mimmi} again.
\qed

\begin{remark} \label{cito2}
\rm The argument of the proof above   can be used to show a similar result on the cones $\U_P$ for a map $f$ conformal in the  product of a standard
sphere  $\mathbb{S}^k$ with an Euclidean space $\R^m$, $k\ge 2$, $m+k\ge 4$. It turns out that $\U_P= T_P\mathbb{S}^k\cup T_P\R^m$. Moreover, all 
arguments of the following Subsection \ref{2952} reduce to a few lines and it is easy to see that a conformal map on a connected open set $\Omega\subset \mathbb{S}^k \times \R^m$ must be the restriction of a local isometry.
\end{remark}

\subsection{The conformal factor $u$} \label{2952}
Here we find all functions $u$ which can be  conformal factors of some conformal maps. We begin by proving 
in the following easy lemma
  that the function $u$ must     be a product. Write $h=(r,y)$
and denote
 by $(h,s)$ points in $\H\times \S$.
\begin{lemma}\label{sssss}
Let $\Omega \subset \H\times\S$ be a connected open  set.
Let $f:\Omega\to f(\Omega)\subset\H\times\S $ be a conformal diffeomorphism with respect to the warped metric $g=g_\H+w^2g_\S$.  Assume that
   $f$ is a product map    of the form either 
\begin{equation}
  \label{1251}
  (h,s) \mapsto (\wt h(h), \wt s(s))  ,\end{equation}
  or
\begin{equation}  \label{1252}
   (h,s)\mapsto (\wt h(s), \wt s(h)),\qquad
    \end{equation} 
    for all $ (h,s)\in \Omega$. 
Then the conformal factor $u$ is a product: $u(h,s)=A(h)B(s)$.
\end{lemma}
\noindent\it Proof. \rm The Cauchy--Riemann system $g(f_*X, f_*X) = u^{-2}g(X, X)$ for all vector field  $X  $ holds.
In case \eqref{1251}, fix  a (lifted) horizontal vector  field $X$.
Observe that $g(X, X)= g_\H(X, X)$ depends on $h$ only. Moreover, by
\eqref{1251},  we have $f_*(X)= \wt h_*(X)$. Therefore
$
g(f_*(X), f_*(X))= g_\H(\wt h_*(X), \wt h_*(X))$ is  a function of $h$ only.
Therefore, by the Cauchy--Riemann system, $u$ depends on $h$ only.

In case \eqref{1252}, which may happen only if $\H$ and $\S$ have the same dimension, given the same $X$ as above, we have $f_*(X) = \wt s_*(X)$, a vertical vector field. Therefore,
by the warped form of $g$,
\[
g(f_*X, f_*X)= g(\wt s_*(X), \wt s_*(X)) = w(\wt h(s))^2 g_\S(\wt s_*(X), \wt s_*(X))=\phi(s)\psi(h),
\]
a product of suitable functions $\phi $ and $\psi$ of $s$ and $h$, respectively.  Therefore the Cauchy--Riemann system gives $u(h,s)=A(h)B(s)$. \endp

\bigskip
Before discussing the system \eqref{pazziata3}, we prove the following proposition.
\begin{proposition}
 \label{www}Let 
$f:\Omega\to f(\Omega)\subset \H \times\S$ be a  conformal
diffeomorphism in the metric $g=g_\H + w^2 g_\S$,   $p\ge
3$. 
Assume that \eqref{caso1} holds. Then $f$ preserves the Ricci tensor. 
\end{proposition}
\noindent\it Proof. \rm We need to prove that  $\Ric(f_*U, f_*V)= \Ric(U, V)$ for all vectors $U, V$. Assumption \eqref{caso1} and the form \eqref{riciao} of the Ricci tensor show that it suffices to assume $U,V\in TS$. In this case we have
\[
\begin{aligned}
 \Ric(f_*U, f_*V)-\Ric(U, V) &=-\a(\a+2)(p-2)(\a+1)^{-2}\big(\wt r^{-2}g 
(f_*U, f_*V)-r^{-2}g(U, V) \big) 
\\& =-\a(\a+2)(p-2)(\a+1)^{-2}\big(\wt r^{-2}u^{-2} -r^{-2} \big)
g(U, V) .
\end{aligned}
\]
To prove the proposition it suffices to show that
$u^{-2}\wt r^{-2}-r^{-2}=0$.
We use the Weyl
 tensor.
 Let $X,Y \in T\S$ be orthogonal vectors. Then, by \eqref{wwi},  \eqref{mimmi} and 
\eqref{riciao} it is easy to see that
\begin{equation}\label{dopps}
 W(X, Y, X, Y)= C_0 r^{-2}g(X, X)  g(Y, Y) ,
\end{equation}
where   $C_0=\displaystyle{ -\frac{q^2+q}{(n-1)(n-2)}\frac{\a(\a+2)}{(\a+1)^2} }<0$, if $\a>0$.
Conformal invariance of $W$ gives
\begin{equation*}
 W(f_*X, f_*Y,f_*X, f_*Y)=u^{-2} W(X, Y, X, Y). 
\end{equation*}
Using \eqref{dopps} in both sides together with the CR system $g(f_*Z, f_*Z)=u^{-2}g(Z, Z)$, $Z=X,Y $, we conclude 
 that
$
 u^{-2}\wt r^{-2}=r^{-2}.
$
Thus, the proposition in proved. \qed
 
Now we are ready to solve  system \eqref{pazziata3}.

\begin{theorem}\label{asr} Let $p\ge 3$ and let   $\Omega\subset\M_0$ be a connected  open set.
Let $f: \Omega\to f(\Omega)\subset \M_0$ be a  smooth diffeomorphism, conformal   in the metric $g$.
Then,  either its conformal factor $u$ is  constant, or it has the form
\begin{equation}\label{tutu}
   u =  a \left( r^2 +\left|
     y-b\right|^2\right)=a \left( |x|^{2(\a+1)} +\left|
     y-b\right|^2\right) = a \|(x, y-b) \|^{2(\a+1)},
\end{equation}
for suitable
$a>0$, $b\in \R^q$, $(x,y)\in\Omega$.
\end{theorem}
\noindent\it Proof. \rm 
Write the map in the  form
$(r,y,\theta)\mapsto (\wt r, \wt y, \wt \theta)$.
We first write    system \eqref{pazziata3} in both cases \eqref{caso1} and \eqref{caso2}.
 Assume that  \eqref{caso1} holds. Then by Proposition \ref{www}, system 
\eqref{pazziata3} becomes  
\begin{equation}\label{pazziatella} 
 (n-2)u^{-1}\D^2  u(U, V) - g(U, V)u^{-2}\big\{ (n-1)|\D u|^2 -
u\Delta u \big\}=0 ,
\end{equation}for any $U,V\in T\Omega$. Here  $\wt r = \wt r(r,y)$.
In case \eqref{caso2} it turns out
 that, given $U=U_\H + U_\S$, we have $(f_* U)_\S = f_*(U_\H)$. Then,
\begin{equation}\label{riso} \begin{split}
 \a(\a+2)&(p-2) (\a+1 )^{-2}
\Big\{r^{-2}g(U_\S, V_\S) - \wt r^{-2 }u^{-2} g\big(U_\H, V_\H\big)\Big\}
\\& =
 (n-2)u^{-1}\D^2  u(U, V) - g(U, V)u^{-2}\big\{ (n-1)|\D u|^2 -
u\Delta u \big\},
\end{split}
\end{equation}
with $\wt r = \wt r(\theta)$.

  Next we start to analyze the systems just obtained.
    The first part of the discussion  is the same for
     case \eqref{pazziatella} and \eqref{riso}.
Indeed, since the connection is Euclidean in variables $r, y_\l$,
in both cases  we have
$
 \p_r\p_\l u=      \D^2 u (\p_r, \p_{\l}) = 0,$  $\l=1,...,q$. We also have
       $\p_\l\p_\m u=\D^2 u(\p_{\l}, \p_{\m}) =0 $,  for all $\l\neq \m $.
Then  $u(\theta,r,y) = F(\theta,r)+\sum_\l G^{(\l)}(\theta,y_\l)$, for suitable
 functions $F, G^{(\l)}$.
 Moreover,  since  $\D_{\p_r}\p_r=0$, $\D_{\p_{\l}}\p_{\l}=0$, both \eqref{pazziatella} and \eqref{riso} give
\begin{equation}\label{dse}
    u_{rr} 
=      \D^2 u (\p_r,\p_r) = \D^2 u (\p_{\l}, \p_{\l} ) = \frac{\p^2 u}{\p y_\l^2},\qquad \l=1,...,q.
\end{equation}
Recall also that, by Lemma \ref{sssss},  $u $ must be   a product.
Thus its form is
\begin{equation}
 \label{gigio}
    u(r,\theta,y) = H(\theta)\Big\{  \frac 1 2 (r^2+|y|^2) +l r
    +\langle m, y\rangle + n\Big\},
\end{equation}
$l,n\in \R$, $m\in \R^p$. Here we used the fact that $\Omega $ is connected.

Next we use condition
$
       \D^2 u (\p_r,X) = 0$, for any $ X$ 
   on the sphere, which holds  in both cases \eqref{pazziatella} and \eqref{riso}.
Let  $X$ be (the lifting of) a vector field on the sphere.
By \eqref{fiorfiore} we get
$
    \p_r X u =\frac 1 r X u 
 $, which gives $    X u (\theta,r,y) = K(\theta,y) r,
$
where   $K$ is a function depending on the vector field  $X$.
Applying   $X$  to   \eqref{gigio} and equating homogeneous powers of $r$,
 we deduce
$
   XH  =0$. Thus $H$ is constant,
    $u$ is constant on the sphere and has the form
\begin{equation}
 \label{gigetta}
    u(r,\theta,y) = \frac 1 2 H(r^2+|y|^2) +L r
    +\langle M, y\rangle+ N,
\end{equation}
for some $M\in \R^q$, $L, N\in \R$.

Next we are ready to rule out case \eqref{riso}.
Indeed,
 letting $U=V=\p_r$ in 
  \eqref{riso}, we get 
\[
 -\a(\a+2)(p-2)(\a+1)^{-2} \wt r^{-2}u^{-2} =
 (n-2)u^{-1} H  -  u^{-2} \big\{ (n-1)|\D u|^2 -
u\Delta u \big\}.
\]
Multiplying by $u^2$ and using the fact that   $\wt r = \wt r(\theta)$ (compare \eqref{zione}),
 we
 get an equation of
 the form
  $\wt r (\theta)^{-2} = \phi(r,y)$, where  $\phi $ is a suitable function.
Therefore it must be $\wt r = $ constant. But this is impossible, because in this case  the map $f$
 would become
singular.

We are left with the study of
 \eqref{pazziatella}. 
Given any unit vector $X\in T\S$ we have
\[
 \D^2 u(X, X)=\D^2 u(\p_r,\p_r).
\]
Taking  the form \eqref{gigetta} of $u$ and \eqref{sforo} into account, we get  $\D^2u(X, X)= \Big(H+\frac{L}{r}\Big)g(X, X)=H+\frac{L}{r}$. Moreover, since   $\D^2 u(\p_r,\p_r)=H$, we conclude that $L=0$. Thus
\begin{equation}
 \label{gigettina}
    u(r,\theta,y) = \frac 1 2 H(r^2+|y|^2)  
    +\langle M, y\rangle+ N,
\end{equation}
$M\in \R^q$, $L, N\in \R$.

Taking the trace of \eqref{pazziatella}, we obtain
\begin{equation}\label{traccia}
 2u\Delta u-n|\D u|^2=0.
\end{equation}
Some  computations based on \eqref{gigettina} and \eqref{tato} give
\begin{equation*}
\begin{split}
     |\D u|^2 & =   |\p_r u|^2 +|\D_y u|^2\\
       & =
     H^2r^2  +|Hy+M |^2
 =
       2 H u - 2N H  +|M|^2, 
 \\   \Delta u  & =  u_{rr}+\frac{p-1}{r} u_r +\Delta_y u
    =   (p-1) H 
   + (q+1) H=nH.
\end{split}
\end{equation*}
Inserting these information into \eqref{traccia}, we easily see that
$
 |M|^2= 2NH.
$
Ultimately,  if   $H=0$, then   $M=0$ and $u =N>0$. If instead 
$H>0$, then we can write
$
     u = \frac{H}{2} \left( r^2 +\left|
     y+\frac{M}{H}\right|^2\right),
$
as desired. \endp

  \subsection{Conclusion of the argument. }\label{triotre}
  Let  $\Omega\subset \H\times_w \S$ be a connected open set.
  Let   $f:\Omega\to f(\Omega)\subset \H\times_w\S$
   be a conformal  diffeomorphism with respect to $g$.  Then, either its conformal
    factor is  constant or it has the form  given in \eqref{tutu}. 
  Recall that the map
  $\Phi(z)= \d_{t\| (x,y-b)  \|^{-2}} (x,y-b)$ has conformal
   factor $u_\Phi (z) = t^{-{\a+1}}\|(x, y-b)\|^{2(\a+1)}$ 
(see Subsection \ref{lasec}, especially equation \eqref{ciccia2}). Write
  $f(z) = F(\Phi(z))$  
and note  that $u_{F\circ\Phi}(z)= u_F(\Phi(z)) u_\Phi(z)$. Then, letting  $t^{-(\a+1)}=a$, the map  $F$ turns out to be a local isometry. The proof is easily concluded, because local isometries are  classified in  Proposition \ref{ioio}.

\section{Umbilical surfaces }\label{ooon}
 \setcounter{equation}{0}
 In this section we prove Theorem \ref{ombelo}.
Let $\Sigma\subset (\M_0,g)$
 be a smooth orientable connected hypersurface. Fix a unit   normal
 vector field $N$.
 Recall that $\Sigma$   is umbilical if at any point $P\in \Sigma$ there is
  $\k(P)\in \R$
  such 
that  the shape operator $L$ satisfies $L(X):=-\D_X N = \k(P) X$, for all $X\in T_P\Sigma$.

Let $p\ge 3$. As discussed in the introduction, the
identification of the cones $\U_P$ and Codazzi equations
 give the following obstruction.
 If  $\Sigma $ is an umbilical surface,
   $P\in \Sigma$ and $N$  is a  normal vector  to $\Sigma$ at $P$, then it must be $N\in
   \U_P$, which means
\begin{equation}\label{casetta}
|N_\H|\,|N_\S|=0.
\end{equation}
Hence, if \eqref{casetta} is not satisfied for a given $N\in T_P\M_0$,
 then  there is no umbilical surface containing $P$ and with normal $N$ at $P$.

Before proving Theorem \ref{ombelo} observe the following facts:

 \begin{remark} \label{ombrello2} \rm  
(1) Since for $p=2$, $q\ge 1$ the manifold $(\M_0, g)$  is flat,
 then for any point $P$ and  $N\in T_P\M_0$ there is    $\Sigma$ umbilical and with normal  $N$ at $P$.

 (2)  The notion of umbilical surface is conformally invariant, while curvature 
depends on the metric.
 The choice of the particular metric $g$ makes all umbilical surfaces to have 
constant curvature. More precisely,  spheres  A1  (as defined in the statement of Theorem \ref{ombelo}) have curvature $1/c$, while planes  A2 and B are geodesic surfaces.

(3) Surfaces A1   can be conformally mapped in surfaces of type
A2, while surfaces of type B cannot.

 (4)  The homogeneous spheres A1 have the same level sets of the
function $\Gamma(z)= \|z\|^{-Q+2}$, $Q= p+(\a+1)q$, which is a singular solution of the equation $\Delta_\a \Gamma=0$  (see \eqref{numero2}) and plays an important role in analysis and potential theory   (see \cite{MM}).
\end{remark}

\noindent\it Proof of Theorem  \ref{ombelo}. \rm The proof will be accomplished in three steps:

 \it Step 1. \rm
Surfaces A1, A2 and B are umbilical.

 \it Step 2.  \rm If $\Sigma $ is umbilical  and
 has normal $\bar N\in T_{\bar P}\H $ at some  $\bar P\in \Sigma$, then $\Sigma$ 
is contained in a surface of type A1 or A2.

 \it Step 3.  \rm
If $\Sigma $ is umbilical  and  has normal $\bar N\in T_{\bar P}\S $ at some  $\bar P\in \Sigma$, then $\Sigma$ is contained in a plane of type B.

\medskip\noindent\it Proof of Step 1.  \rm We start from type A1. Without loss of generality we take $c=1$ and $b=0$. so that our surface $\Sigma$ has equation $|x|^{2(\a+1)}+|y|^2=1$
 In the warped model with metric \eqref{harpo4},
  a  unit normal vector field  has the form  
$
N= -  (r\p_r + y^\l \p_{ \l}).
$
Let $P\in\Sigma$ and $U\in T_P\Sigma$.  By linearity of the shape operator, it suffices to consider separately the cases 
 $U\in T_P \S$ and $U=
a\p_r +c_\l\p_\l$, where $ar + \sum c_\l y^\l=0$.
 In the first case,
$
L(U) = -\D_U N=   \D_U (r\p_r + y^\l \p_{\l})=U,
$
in view of \eqref{fiorfiore}. In the second case, if   $U= a\p_r +c_\l\p_\l$,
then $
L(U) =  -\D_UN=\D_{ a\p_r +c_\l\p_\l}(r\p_r + y^\m\p_\m)=
U,
$
because in these variables the connection is Euclidean.
The proof for planes  A2 is analogous and we omit it.

Next we pass to Type  B.   Here $\Sigma $ has equation $\sum_{k}a_k x_k=0$.  
Assume that $\sum
a_k^2 = 1$.  We use Cartesian coordinates $(x^j, y^\l)$ and the metric $g$.
A unit normal  vector field  to $\Sigma$ is  
 $N = (\a+1)^{-1}|x|^{-\a} a_k\p_{x_k}.
 $
 Again by linearity of $L$, it suffices to consider separately 
 vectors of the form  $U= \p_{\l}$, 
with  $\la=1, \dots, q$, 
 and  $U = U^j \p_{j}$, which are tangent provided $U^k a_k=0$.
 In first case,
 $
 L(\p_{\l}) =  -(\a+1)^{-1} \D_{\p_{\l}}
 |x|^{-\a} a_k\p_{k}=0,
 $
because $\D_{\p_{\l}}\p_{k}=0$. In the second case, since the
Christoffel symbols of the metric   $(\a+1)^2|x|^{2\a}|dx|^2$ in
$\R^p\setminus\{0\}$
 are $
\Gamma_{ij}^k   = \a|x|^{-2}\big\{\d_{ik }x_j +\d_{jk}x_i -\d_{ij}
x_k \big\}$, we get
\[
\begin{split}
L(U^j\p_j) & =  -  U^j\D_{\p_{j}} (\a+1)^{-1} |x|^{-\a}a_k\p_{k}
\\&
=- (\a+1)^{-1} U^j a_k \big\{
 -\a|x|^{-\a-2} x_j\p_k + |x|^{-\a}\Gamma_{jk}^i\p_i\big\}
\\&
   =  \a(\a+1)^{-1}|x|^{-\a-2}\{- (a_k x_k) U + U^k a_k x_i\p_i  \} = 0,
\end{split}
\]
because $\langle a, x\rangle=0$ and   $U$ is tangent to the plane.
Thus $L(U)=0$. We have proved that $\Sigma$ is a geodesic surface.
\qed

\medskip\noindent\it Proof of Step 2. \rm  Let $\bar P\in \M_0$ and let
 $\bar N = \bar N_\H\in T_{\bar P} \H$.
 Let  $\Sigma$ be an umbilical surface with normal $\bar N  $
  at $\bar  P$ with respect to $g$.
 Examples of
 surfaces  A1 and A2 show that there 
is at least one surface with these properties.
 We want 
to show that $\Sigma$ is contained in a surface of type A1 or A2. Denote by $N$ the
 unit normal vector field  to $\Sigma$ which agrees with $\bar N$ at $\bar P.$
  By continuity and by \eqref{casetta} it must be that
$N= N_\H\in T_P\H$
 at any point   $P\in \Sigma$.

We prove first that  $\Sigma$ has constant curvature  $\k$.
 It suffices to prove that
$U\k=0$ for any  vector $U$ tangent to $\Sigma$.
 Since $N= N_\H$, \eqref{mimmi} gives  $R(N, U, V, W) =0
$, for any $U, V, W$ orthogonal to $N$.  Thus, Codazzi equations \eqref{paippo} show that $\k$ must be constant.

In the warped model $(r,y,\theta)$,  the  vector field $N$ has the form
  $
 N= a\p_r + N^\l\p_{\l},$  where $a, N^\l$  are 
  suitable functions on $\Sigma$.
Since $\Sigma $ is umbilical, given any  $V= V_\H +V_\S \in T_P\Sigma$,
 it must be that
 \begin{equation}
 \label{cicao}
 \begin{split}
 -\k V_H-\k V_\S =-\k V =\D_V N & =  (Va) \p_r + a\D_V\p_r +  (VN^\l)\p_\l + N^\l \D_V \p_{\l}
 \\&=  (Va)\p_r +\frac ar V_\S +( VN^\l)\p_{\l} ,
 \end{split}\end{equation}
 where we used \eqref{cane} and \eqref{fiorfiore}.
 Comparing like terms, we get
\begin{equation}
\label{cisa}
-\k = a/r ,\quad \Rightarrow
N= -r\k \p_r + N^\l \p_{\l}.
\end{equation}
Since  $|N|=1$, it must be that
\begin{equation}
\label{novo2}
r^2 \k^2 + \sum_\l (N^\l)^2=1.
\end{equation}
Write  in  \eqref{cicao} $V_\H = V^\l \p_{\l} + V^r\p_r$ and take
components along    $\p_{\l}$. Thus
 \begin{equation}
 \label{tutto2}
 -\k V^\l = VN^\l, \quad\forall\; V\in T\Sigma.
 \end{equation}
There are two cases:
if  $\k=0$, then  \eqref{tutto2} implies that $N^\l $ is constant. Therefore $\Sigma $ is contained in
 the plane of equation $N^\l y^\l=$constant.
If instead  $\k\neq 0$,   \eqref{tutto2} gives
$
V(y^\l + \k^{-1}N^\l)=0$. Therefore $ \quad y^\l + \k^{-1}N^\l=b_\l,
$ where $b_\l$   is  a constant. Thus   \eqref{novo2} becomes
 $ 1= r^2\k^2 + \k^2
\sum_\l (b_\l- y_\l)^2, $ as desired.
\qed

\smallskip \noindent\it Proof of Step 3. \rm 
Let $P\in \M_0$ and let   $\bar N =\bar  N_\S\in T_{\bar  P}\S $ be a
unit vector. Consider an umbilical surface $\Sigma$   with normal
$\bar N$ at $\bar  P$. Surfaces of type B show that there is at least one
surface with this property.
 Our aim is to show that $\Sigma$ is
contained in a plane of type B.

Let  $N$ be the unit
normal  to $\Sigma$ which agrees with   $\bar  N $  at $\bar P=(\bar r, \bar  y, \bar  \theta)$.
Since   $\Sigma $ is umbilical, by  \eqref{casetta} and by continuity, it must be that  $N_\H=0$ for all $P= (r,y,\theta)\in \Sigma$.
Therefore, given a local frame     $X_j$, $j=1, \dots, p-1$,
on the sphere     $\mathbb{S}^{p-1}$ around
 $\bar \theta  $,   $N$ has the form
 $
 N=   \sum_{j=1}^{p-1} b_j(r,y,\theta) X_{{j}}.
 $

 Next take the tangent vector  $\p_{y_1}\in T_{ P}\Sigma ,$ for any  $P $ close to $\bar P$.
   Let  $\k$  be the curvature of  $\Sigma$. Then
  \[
  \k \p_{y_1}  = L(\p_{y_1})
  =- \sum_{j\le p-1}\big\{(\p_{y_1}b_j)
  X_j
  +b_j \D_{\p_{y_1}}X_j \big\}= -\sum_{j\le p-1} (\p_{y_1}b_j)
  X_j
  ,
  \]
  by \eqref{fiorfiore}.
Therefore it must be that   $\k=0$. Hence $\Sigma$ is a  geodesic
surface. Thus it  must be contained  in the plane
of equation  $\sum_{k}\bar  N^k x_k=0$, which is by the previous Step 1 a geodesic surface too.
\qed

\section*{Appendix}
\setcounter{equation}{0}\renewcommand{\theequation}{A.\arabic{equation}}
 We collect here some standard formulas on warped products. See \cite[Chapter 7]{ON} for a complete discussion.
Let  $(\H, g_\H)$ and
$(\S,g_\S) $ be Riemannian manifolds. Given
$w:\H\to \left]0,+\infty\right[$,  the warped product $\H\times_w \S $ is the
manifold $\H\times\S$ equipped with the metric $g= g_\H + w^2 g_S$.
Given any $P= (h,s)$,
Decompose as usual $T_{P}M$ as the orthogonal sum of $T_P \H$ and $T_P\S$,
  the   lifts  at $P$ of $T_h\H$
and $T_s \S$, respectively.
We use the same notation for a  vector and its lifting.
Lifting of vector fields on $\H$ and on $\S$ are usually denoted by $\L(\H)$ and $\L(\S)$. They are often called
\emph{lifted horizontal} or \emph{lifted vertical} vector fields.  Vector fields ant their liftings are denoted by the same symbol. Observe that
for a function $\phi$ depending on  $h$ only, the gradient  $\grad \phi$ of $\phi$ in the metric $g$
 is nothing but the obvious
lifting of $\grad_{g_\H}\phi$.

Next, let $\D^\H$, $\D^\S$ and $\D$ be the Levi Civita connections
on $\H$, $\S$ and $\H \times_w\S$, respectively.  Then,
the following formulas hold  (below $A, B, C, D\in \L(\H)$ and $X, Y, Z, V\in \L(\S)$).
\begin{equation}
\begin{split}
\label{alfa}
 \D_{A}B & =\D^\H_{A}B,
\\ \D_{A}X & =\D_{X}{A}=  w^{-1}(A w) X,
\\ \D_{X}Y &= \D^\S_{X}Y -g(X, Y)w^{-1}\grad w.
\end{split}
\end{equation}
Therefore, given $u:\H\times\S\to \R$, with slight abuse of notation,
\begin{equation}
\begin{split}
\label{beta} \D^2 u(A, B) &=
  (\D^\H)^2u (A, B),
\\
\D^2 u(X, Y) &=
  (\D^\S)^2u (X, Y)+w^{-1}  g(X, Y) (\grad w)u  ,
\\
\Delta u &= \Delta_\H u + w^{-2}\Delta_\S u + \mathrm{dim}(\S)w^{-1}
(\grad w) u.
\end{split}
\end{equation}

A computation using \eqref{alfa} provides also
\begin{equation*}
\begin{split}
R(A, B )C & = R_\H(A, B)C,
\\
R(A , X) B   & =  w^{-1} \D^2 w (A, B)X,
\\  R(X, Y) Z & =R_\S(X, Y)Z +
w^{-2}|\D w|^2\big\{g(X, Z)Y - g(Y, Z)X\big\}.
\end{split}
\end{equation*}
Therefore,
\begin{equation}
\begin{split}
\label{beta3} R(D, C, A, B) & :=g(D, R(A, B) C) =R_\H(D, C, A, B),
\\
R( Y,B,  A, X) & =
 w^{-1}\D^2 w (A, B) g(X,Y)
\\
R(V, Z, X, Y)& = w^2R_\S (V, Z, X, Y)
\\
\qquad & + w^{-2}|\D w |^2\big\{ g(X, Z)g(Y, V)- g(Y,
Z)g(X, V) \big\}.
\end{split}
\end{equation}
The remaining nonzero components of $R$ can be obtained by the
standard symmetries  $R_{abcd}= -R_{bacd}= R_{cdab}$ of the curvature tensor
 $R$.

\small

\medskip
 \noindent Daniele Morbidelli: Dipartimento di Matematica,\\
 Universit\`a di Bologna.\\
 Piazza di Porta San Donato, 5,\\
 40127  Bologna, Italy.\\
E mail: \tt morbidel@dm.unibo.it

%
\end{document}